\documentclass [12pt]{article}
\usepackage{graphicx,amssymb,amsfonts,latexsym,amsmath,amsthm,times}
\usepackage{epsfig}
\usepackage{color}
\usepackage[english]{babel}
\usepackage[a4paper]{geometry}
\def\lanbox{\hbox{$\, \vrule height 0.25cm width 0.25cm depth 0.01cm \,$}}

\numberwithin{equation}{section}

\begin{document}

\vspace*{1.4cm}

\normalsize \centerline{\Large \bf EXISTENCE OF
STATIONARY SOLUTIONS FOR  SOME}

\medskip

\normalsize \centerline{\Large \bf 
SYSTEMS OF INTEGRO-DIFFERENTIAL EQUATIONS}

\medskip

\centerline{\Large\bf WITH LAPLACE AND BI-LAPLACE OPERATORS}

\vspace*{1cm}

\centerline{\bf  Vitali Vougalter$^{1*}$,  Vitaly Volpert$^{2, 3}$   }

\vspace*{0.5cm}

\centerline{$^{1}$ Department of Mathematics, University
of Toronto}

\centerline{Toronto, Ontario, M5S 2E4, Canada}

\centerline{ e-mail: vitali@math.toronto.edu}

\centerline{$^*$ Corresponding author}

\medskip

\centerline{$^2$ Institute Camille Jordan, UMR 5208 CNRS,
University Lyon 1}

\centerline{ Villeurbanne, 69622, France}

\centerline{$^3$ Peoples' Friendship University of Russia, 6 Miklukho-Maklaya
St,}

\centerline{Moscow, 117198, Russia}

\centerline{e-mail: volpert@math.univ-lyon1.fr}

\medskip


\vspace*{0.25cm}

\noindent {\bf Abstract:}
The article is devoted to the solvability of a system of
integro-differential equations in the case of the difference of the standard Laplacian and 
the bi-Laplacian in the diffusion terms.
The proof of the existence of solutions is based
on a fixed point technique. We use the solvability conditions for the
elliptic operators without the Fredholm property in unbounded domains.

\vspace*{0.25cm}



\vspace*{0.5cm}

\bigskip

\bigskip


\setcounter{section}{1}

\centerline{\bf 1. Introduction}

\medskip

In the present work we address the existence of stationary
solutions of the following system of $N\geq 2$ integro-differential equations in
${\mathbb R}^{d}, \ 5\leq d\leq 7$
\begin{equation}
\label{h}
\frac{\partial u_{m}}{\partial t}=
D_{m}[\Delta-{\Delta}^{2}]u_{m}+
\int_{{\mathbb R}^{d}}K_{m}(x-y)g_{m}(u(y,t))dy + f_{m}(x),
\end{equation}
where $1\leq m\leq N$, which is relevant to the cell population dynamics.
The results of the article are derived in this particular range of the
values of the dimension, which is based on the
solvability of the linear Poisson type equations (\ref{lpm}) and the
applicability of the Sobolev embedding (\ref{e}).
The solvability of the single equation analogous to
(\ref{h}) in $H^{3}({\mathbb R}^{5})$ involving a single Laplacian in the diffusion term
was studied in ~\cite{VV15} (see also ~\cite{VV111}).
Note that the space variable $x$ in our system is
correspondent to the cell genotype, the functions
$u_{m}(x,t)$ describe the cell density distributions for various groups of
cells as functions of their genotype and time,
$$
u(x,t)=(u_{1}(x,t), u_{2}(x,t),...,u_{N}(x,t))^{T}.
$$
The right side of the system of equations (\ref{h}) describes the evolution of
cell densities by virtue of the cell proliferation, mutations and the cell influx
or efflux. The diffusion terms with positive
coefficients $D_m$ contain  the difference of the Laplace operator and the bi-Laplacian. 
The integral production terms describe large mutations.
The functions $g_{m}(u)$ stand for the rates of the cell birth depending on $u$
(density dependent proliferation), and the kernels $K_{m}(x-y)$ denote
the proportions of the newly born cells, which  change their genotypes from $y$ to $x$.
Assume that they depend on the distance between the genotypes.
The functions $f_{m}(x)$ designate the influxes or effluxes of cells for
different genotypes.

The bi-Laplacian in the diffusion terms is important to describe the long range interactions in biological systems (see ~\cite{M03}). 
The diffusion operator involving  the difference of the Laplace and bi-Laplace operators can appear in several modelling contexts. 
In the evolutionary dynamics,  the Laplacian term describes the small random mutations in the genotype space, 
while the higher-order bi-Laplacian may represent the long-range smoothing effects or higher-order mutation processes. The operators
of this kind  appear in the ecological models where individuals disperse through a combination of the local diffusion and longer-range redistribution.  From another point of view,  the operator $\Delta-{\Delta}^{2}$ can be considered  as a local approximation
 of nonlocal dispersal operators  derived by the Taylor expansion of convolution kernels. Note that the higher-order diffusion terms
arise frequently in the pattern formation theory and in models with higher-gradient regularization, where they reflect additional spatial interactions or energetic penalties associated with the curvature of the density distribution. The global well-posedness of an
integro-differential equation with the bi-Laplacian and transport was established in ~\cite{EV25}.

Let us set here for the simplicity all $D_{m}=1$. This can be accomplished by the division of each $m$th equation of our system of equations by $D_{m}$. We establish the existence of solutions of the problem
\begin{equation}
\label{p}
[\Delta-{\Delta}^{2}]u_{m}+\int_{{\mathbb R}^{d}} K_{m}(x-y)g_{m}
(u(y))dy +f_{m}(x) = 0, \quad 5\leq d\leq 7,
\end{equation}
where $1\leq m\leq N$.
We consider our system in the whole ${\mathbb R}^{d}$. This enables us to exclude the influence of the boundary conditions.
The corresponding operator studied in the whole space is non-Fredholm. This absence of
the Fredholm property naturally appears in the mathematical formulation relevant to our biological
problem. Consequently,  the traditional approaches
of the nonlinear analysis may not be applied.  Our argument relies on
the solvability conditions for the operator which fails to satisfy the Fredholm property and
the method of contraction mappings.

Let us consider the equation
\begin{equation}
\label{eq1}
 -\Delta u + V(x) u - a u=f,
\end{equation}
where $u \in E= H^{2}({\mathbb R}^{d})$ and  $f \in F=
L^{2}({\mathbb R}^{d}), \ d\in {\mathbb N}$, $a$ is a constant and
the scalar potential function $V(x)$ either vanishes in the whole space
or tends to $0$ as $|x|\to \infty$. Such model problem is discussed here
to illustrate the particular features of the equations involving the non-Fredholm
operators, the techniques used to solve them and the preceding results.
If $a \geq 0$, the origin belongs to the essential spectrum of the
operator $A : E \to F$, which corresponds to the left side of 
(\ref{eq1}). Consequently, this operator does not
satisfy the Fredholm property. Its image is not closed, for $d>1$
the dimension of its kernel and the codimension of its image are
not finite. Our article is devoted to the studies of the certain properties
of the operators of this kind. The elliptic problems containing the non-Fredholm
operators were considered actively in recent years.
Approaches in weighted Sobolev and H\"older spaces were developed in
~\cite{Amrouche1997}, ~\cite{Amrouche2008}, ~\cite{Bolley1993},
~\cite{Bolley2001}, ~\cite{B88}. The Schr\"odinger type operators, which do not satisfy the
Fredholm property were studied with the methods of the spectral and the
scattering theory in  ~\cite{EV21}, ~\cite{EV22}, ~\cite{V2011}, ~\cite{VV08}, ~\cite{VV14}, ~\cite{VV22}.
The nonlinear non-Fredholm elliptic equations were covered in ~\cite{EV22}, ~\cite{VV111}, ~\cite{VV14}, ~\cite{VV15}, ~\cite{VV21}.
Work ~\cite{DH22} is devoted to a two phase boundary obstacle-type problem for the bi-Laplacian. 
Article ~\cite{DKV19} deals with  the limit behaviour of a singular perturbation problem for the biharmonic operator. 
The important applications to the theory of reaction-diffusion
type equations were investigated in ~\cite{DMV05}, ~\cite{DMV08}. Fredholm
structures, topological invariants and applications were discussed in
~\cite{E09}. Works ~\cite{GS05} and ~\cite{RS01} are significant for the
understanding of the Fredholm and properness properties of the quasilinear
elliptic systems of the second order and of the operators of this kind
on ${\mathbb R}^{N}$.
The non-Fredholm operators arise also when studying the wave systems with
an infinite number of localized traveling waves (see ~\cite{AMP14}). Standing lattice solitons in the discrete NLS equation 
with saturation were covered in ~\cite{AKLP19}.
Particularly, when $a$ is trivial, the operator $A$ is Fredholm in certain properly chosen
weighted spaces (see \cite{Amrouche1997}, \cite{Amrouche2008},
\cite{Bolley1993}, \cite{Bolley2001}, \cite{B88}). However, the situation when
$a \neq 0$ is considerably
different and the methods developed in these articles cannot be used.

\medskip

We set $K_{m}(x) = \varepsilon_{m} H_{m}(x)$, where
$\varepsilon_{m} \geq 0$, so that
\begin{equation}
\label{emax}
\varepsilon:=\hbox{max}_{1\leq m\leq N}\varepsilon_{m}.
\end{equation}

\bigskip

\noindent
{\bf Assumption 1.1.}  {\it Let
$1\leq m\leq N$, the functions
$f_{m}(x): {\mathbb R}^{d}\to {\mathbb R}, \ 5\leq d\leq 7$ are nontrivial for a
certain $m$, such that
$$
f_{m}(x)\in L^{1}({\mathbb R}^{d})\cap L^{2}({\mathbb R}^{d}).
$$
In addition, we assume that
$H_{m}(x): {\mathbb R}^{d}\to {\mathbb R}$, so that
$$
H_{m}(x)\in L^{1}({\mathbb R}^{d})\cap  L^{2}({\mathbb R}^{d}).
$$
Moreover,
\begin{equation}
\label{h2}
H^{2}:=\sum_{m=1}^{N}\|H_{m}(x)\|_{L^{1}({\mathbb R}^{d})}^{2}>0
\end{equation}
and}
\begin{equation}
\label{q2}
Q^{2}:=\sum_{m=1}^{N}\|H_{m}(x)\|_{L^{2}({\mathbb R}^{d})}^{2}>0.
\end{equation}

\bigskip

In the article we work in the space of dimension $5\leq d\leq 7$. This is related to the
solvability conditions for the linear Poisson type problem (\ref{lp}) formulated in
Lemma 4.1 below and to the applicability of the Sobolev embedding (\ref{e}). From the perspective of the practical applications, 
the space dimensions are not limited to $5\leq d\leq 7$, since the space variable here is correspondent to the cell
genotype but not to the usual physical space.
For the technical purposes, we use the Sobolev space
\begin{equation}
\label{h4}
H^{4}({\mathbb R}^{d}):=\big\{\phi(x):{\mathbb R}^{d}\to {\mathbb R} \ | \
\phi(x)\in L^{2}({\mathbb R}^{d}), \ {\Delta}^{2} \phi(x) \in L^{2}
({\mathbb R}^{d}) \big\}.
\end{equation}
Such space (\ref{h4}) is equipped with the norm
\begin{equation}
\label{n}
\|\phi\|_{H^{4}({\mathbb R}^{d})}^{2}:=\|\phi\|_{L^{2}({\mathbb R}^{d})}^{2}+
\|{\Delta}^{2} \phi \|_{L^{2}({\mathbb R}^{d})}^{2}.
\end{equation}
For a vector function
$$
u(x)=(u_{1}(x), u_{2}(x), ..., u_{N}(x))^{T},
$$
in the present work we will use the norms
\begin{equation}
\label{u1N}
\|u\|_{H^{4}({\mathbb R}^{d}, {\mathbb R}^{N})}^{2}:=\sum_{m=1}^{N}\|u_{m}\|_{H^{4}({\mathbb R}^{d}) }^{2}=
\sum_{m=1}^{N}\{\|u_{m}\|_{L^{2}({\mathbb R}^{d}) }^{2}+\|{\Delta}^{2} u_{m}\|_{L^{2}({\mathbb R}^{d})}^{2}\}
\end{equation}
and
\begin{equation}
\label{l2vn}
\|u\|_{L^{2}({\mathbb R}^{d}, {\mathbb R}^{N})}^{2}:=\sum_{m=1}^{N}
\|u_{m}\|_{L^{2}({\mathbb R}^{d})}^{2}.
\end{equation}
By means of the standard Sobolev embedding  in dimensions $d\leq 7$,
\begin{equation}
\label{e}
\|\phi\|_{L^{\infty}({\mathbb R}^{d})}\leq c_{e}\|\phi\|_{H^{4}({\mathbb R}^{d})}.
\end{equation}
Here $c_{e}>0$ is the constant of the embedding.
When all the nonnegative parameters $\varepsilon_{m}$ are trivial, we arrive at
the linear Poisson type equations with $x\in {\mathbb R}^{d}, \ 5\leq d\leq  7$, namely
\begin{equation}
\label{lpm}
[-\Delta+{\Delta}^{2}]u_{m}(x)=f_{m}(x), \quad 1\leq m\leq N.
\end{equation}
By virtue of Lemma 4.1 further down along with Assumption 1.1, each problem
(\ref{lpm}) admits a unique solution
$$
u_{0, m}(x)\in H^{4}({\mathbb R}^{d}), \quad 1\leq m\leq N, \quad 5\leq d\leq  7.
$$
 Recall the definition of the norm (\ref{u1N}). Hence,
\begin{equation}
\label{u0v}
u_{0}(x):=(u_{0,1}(x), u_{0,2}(x), ..., u_{0,N}(x))^{T}\in
H^{4}({\mathbb R}^{d}, {\mathbb R}^{N}), \quad 5\leq d\leq  7.
\end{equation}
According to Assumption 1.1, the functions $f_{m}(x)$ do not vanish identically in the
whole space for some $1\leq m\leq N$. Therefore, vector function (\ref{u0v}) is nontrivial.
Let us look for the resulting solution of the nonlinear system of equations
(\ref{p}) as
\begin{equation}
\label{r}
u(x)=u_{0}(x)+u_{p}(x),
\end{equation}
where
$$
u_{p}(x):=(u_{p,1}(x), u_{p,2}(x),...,u_{p,N}(x))^{T}.
$$
Then we easily obtain the perturbative system of equations
\begin{equation}
\label{pert}
[-\Delta+{\Delta}^{2}]u_{p,m}(x)=\varepsilon_{m} \int_{{\mathbb R}^{d}}
H_{m}(x-y)g_{m}(u_{0}(y)+u_{p}(y))dy
\end{equation}
with $1\leq m\leq N, \quad 5\leq d\leq 7$.

We introduce a closed ball in the Sobolev space
\begin{equation}
\label{b}
B_{\rho}:=\{u\in H^{4}({\mathbb R}^{d}, {\mathbb R}^{N}) \ | \ \|u\|_
{H^{4}({\mathbb R}^{d}, {\mathbb R}^{N})}\leq \rho \}, \quad 0<\rho\leq 1, \quad 5\leq d\leq 7.
\end{equation}
Let us seek the solution of problem (\ref{pert}) as the fixed point of the
auxiliary nonlinear system
\begin{equation}
\label{aux}
[-\Delta+{\Delta}^{2}]u_{m}(x)=\varepsilon_{m} \int_{{\mathbb R}^{d}}
H_{m}(x-y)g_{m}(u_{0}(y)+v(y))dy,
\end{equation}
where $1\leq m\leq N, \ 5\leq d\leq 7$ in ball (\ref{b}). For a
given vector function $v(y)$ this is a system of equations with respect to
$u(x)$.
The left side of each equation in (\ref{aux}) contains the operator which
does not satisfy the Fredholm property
\begin{equation}
\label{l}
l:=-\Delta+{\Delta}^{2}: H^{4}({\mathbb R}^{d})\to
L^{2}({\mathbb R}^{d}).
\end{equation}
It is the differential operator with the symbol $|p|^{2}+|p|^{4}$,
so that 
$$
l\phi(x)=\frac{1}{(2\pi)^{\frac{d}{2}}}\int_{{\mathbb R}^{d}}
(|p|^{2}+|p|^{4})\widehat{\phi}(p)e^{ipx}dp, \quad
\phi(x)\in H^{4}({\mathbb R}^{d})
$$
with the standard Fourier transform defined in (\ref{f}). Obviously,
the essential spectrum of (\ref{l}) fills the nonnegative semi-axis
$[0, +\infty)$.
Thus, such operator does not have a bounded inverse. The analogous situation
arised in previous articles ~\cite{VV111} and ~\cite{VV14} but as distinct from the
present case, the equations discussed there required the orthogonality
relations.  Persistence of pulses for certain local reaction-diffusion problems via
the fixed point technique was studied in ~\cite{CV21}.
But the Schr\"odinger type operator contained in the nonlinear
equation there possessed the Fredholm property.

Let us introduce the closed ball in the
space of $N$ dimensions as
\begin{equation}
\label{i}
I:=\{z\in {\mathbb R}^{N} \ | \ |z|_{{\mathbb R}^{N}}\leq
c_{e}\|u_{0}\|_{H^{4}({\mathbb R}^{d}, {\mathbb R}^{N})}+c_{e} \}.
\end{equation}
Here and below $|.|_{{\mathbb R}^{N}}$ will denote the length of a vector in
${\mathbb R}^{N}$.
The closed ball $D_{M}$ in the space of $C^{2}(I,{\mathbb R}^{N})$
vector functions is
\begin{equation}
\label{M}
\{g(z):=(g_{1}(z), g_{2}(z),..., g_{N}(z))\in C^{2}(I,{\mathbb R}^{N}) \ | \
\|g\|_{C^{2}(I,{\mathbb R}^{N})}\leq M \},
\end{equation}
where  $M>0$.
In such context the norms
\begin{equation}
\label{gN}
\|g\|_{C^{2}(I,{\mathbb R}^{N})}:=\sum_{m=1}^{N}\|g_{m}\|_{C^{2}(I)},
\end{equation}
\begin{equation}
\label{gn}
\|g_{m}\|_{C^{2}(I)}:=\|g_{m}\|_{C(I)}+\sum_{n=1}^{N}\Big\|\frac{\partial g_{m}}
{\partial z_{n}}\Big\|_{C(I)}+\sum_{n,l=1}^{N}\Big\|\frac{\partial^{2}g_{m}}
{\partial z_{n}\partial z_{l}}\Big\|_{C(I)}
\end{equation}
with $\|g_{m}\|_{C(I)}:=\hbox{max}_{z\in I}|g_{m}(z)|$.
Let us impose the following auxiliary conditions on the nonlinear part of the system
of equations (\ref{p}). From the point of view of the applications in biology,
$g_{m}(z)$ can be, for instance the quadratic functions,  describing the
cell-cell interactions.

\bigskip

\noindent
{\bf Assumption 1.2.} {\it Let $1\leq m\leq N$. Assume that each
$g_{m}: {\mathbb R}^{N}\to {\mathbb R}$ is such that
$g_{m}(0)=0$ and $\nabla g_{m}(0)=0$. Furthermore, $g\in D_{M}$ and
it does not vanish identically in the ball $I$}.

\bigskip

We use the technical Assumptions 1.1 and 1.2 in the proofs of our main
propositions. It is not clear at the moment if there is a more efficient way to
analyze our problem  which would enable us to weaken these
conditions.

Let us introduce the operator $T_g$, such that $u =T_g v$, where $u$ is a
solution of system (\ref{aux}). Our first main statement is
as follows.

\bigskip

\noindent
{\bf Theorem 1.3.} {\it Suppose Assumptions 1.1 and 1.2 hold. Then for every
$\rho\in (0, 1]$ problem
(\ref{aux}) defines the map $T_{g}: B_{\rho}\to B_{\rho}$, which is a strict
contraction for all
$$
0<\varepsilon\leq
\frac{\rho}{M(\|u_{0}\|_{H^{4}({\mathbb R}^{d}, {\mathbb R}^{N})}+1)^{2}}\times
$$
\begin{equation}
\label{eps}
\Bigg[H^{2}(\|u_{0}\|_{H^{4}({\mathbb R}^{d}, {\mathbb R}^{N})}+1)^{\frac{8}{d}-2}\Bigg(\frac{|S^{d}|}{4}\Bigg)^{\frac{4}{d}}
\frac{d}{(d-4)(2\pi)^{4}}
+Q^{2}
\Bigg]^{-\frac{1}{2}}.
\end{equation}
The unique fixed point $u_{p}$ of
such map $T_{g}$ is the only solution of the system of equations (\ref{pert}) in $B_{\rho}$.}

\bigskip

Here and below $S^{d}$ denotes the unit sphere in the space of $5\leq d\leq 7$
dimensions centered at the origin and $|S^{d}|$ stands for its Lebesgue measure.

Note that the constants $\varepsilon, \ H, \ Q$ in formula (\ref{eps}) are defined in
(\ref{emax}),  (\ref{h2}) and (\ref{q2}).

Evidently, the resulting solution $u(x)$ of system (\ref{p})
given by (\ref{r}) will not vanish identically in ${\mathbb R}^{d}$ since
the influx/efflux terms $f_{m}(x)$ are nontrivial for a certain
$1\leq m\leq N$ and all $g_{m}(0)=0$ according to our assumptions. We have
the following auxiliary lemma.

\bigskip

\noindent
{\bf Lemma 1.4.} {\it Suppose $R\in (0, +\infty)$ and $5\leq d\leq 7$. Consider the
function
$$
\varphi(R):=\alpha R^{d-4}+\frac{1}{R^{4}}, \quad \alpha>0.
$$
It achieves its minimal value at \
$\displaystyle{R^{*}:=\Bigg({\frac{4}{\alpha (d-4)}}\Bigg)^{\frac{1}{d}}}$,
which is given by}
$$
\varphi(R^{*})=\Bigg(\frac{\alpha}{4}\Bigg)^{\frac{4}{d}}\frac{d}{(d-4)^{\frac{d-4}{d}}}.
$$

\bigskip

The second main result of the work is devoted to the continuity of the resulting
solution of the system  of equations (\ref{p}) given by formula (\ref{r}) with
respect to the nonlinear vector function $g$. Let us introduce the following
positive technical quantity
$$
\kappa:=
$$
\begin{equation}
\label{kap}
M(\|u_{0}\|_{H^{4}({\mathbb R}^{d}, {\mathbb R}^{N})}+1)
\Bigg\{\frac{H^{2}(\|u_{0}\|
_{H^{4}({\mathbb R}^{d}, {\mathbb R}^{N})}+1)^{\frac{8}{d}-2}d}
{(d-4){(2 \pi)}^{4}}
\Bigg(\frac{|S^{d}|}{4}\Bigg)^{\frac{4}{d}}+
Q^{2}\Bigg\}^{\frac{1}{2}}.
\end{equation}

\bigskip

\noindent
{\bf Theorem 1.5.} {\it Let $j=1,2$, the conditions of Theorem 1.3 are valid,
so that
$u_{p,j}$ is the unique fixed point of the map
$T_{g_{j}}: B_{\rho}\to B_{\rho}$, which is a strict contraction for all the values
of $\varepsilon$, which satisfy inequality (\ref{eps}) and the resulting solution of 
system  (\ref{p}) with $g(z)=g_{j}(z)$ equals to
\begin{equation}
\label{jres}
u_{j}(x):=u_{0}(x)+u_{p, j}(x).
\end{equation}
Then for all the values of the parameter $\varepsilon$ satisfying bound (\ref{eps}),
the estimate
\begin{equation}
\label{cont}
\|u_{1}-u_{2}\|_{H^{4}({\mathbb R}^{d}, {\mathbb R}^{N})}\leq \frac{\varepsilon \kappa}
{M(1-\varepsilon \kappa)}(\|u_{0}\|_{H^{4}({\mathbb R}^{d}, {\mathbb R}^{N})}+1)
\|g_{1}-g_{2}\|_{C^{2}(I, {\mathbb R}^{N})}
\end{equation}
holds.}

\bigskip

We turn our attention to the establishing of the validity our first main statement.

\bigskip


\setcounter{section}{2}
\setcounter{equation}{0}

\centerline{\bf 2. The existence of the perturbed solution}

\bigskip

\noindent
{\it Proof of Theorem 1.3.} We choose an arbitrary  vector function
$v\in B_{\rho}$ and denote the terms contained in the integral expressions
in the right side of the system of equations  (\ref{aux}) as
$$
G_{m}(x):=g_{m}(u_{0}(x)+v(x)), \quad 1\leq m\leq N.
$$
Let us use the standard Fourier transform 
\begin{equation}
\label{f}
\widehat{\phi}(p):=\frac{1}{(2\pi)^{\frac{d}{2}}}\int_{{\mathbb R}^{d}}\phi(x)e^{-ipx}dx,
\quad p\in {\mathbb R}^{d}, \quad 5\leq d\leq 7.
\end{equation}
Evidently, the inequality
\begin{equation}
\label{fub}
\|\widehat{\phi}\|_{L^{\infty}({\mathbb R}^{d})}\leq \frac{1}{(2\pi)^{\frac{d}{2}}}
\|\phi\|_{L^{1}({\mathbb R}^{d})}
\end{equation}
holds.
We apply (\ref{f}) to both sides of problem (\ref{aux}). This gives us
$$
\widehat{u_{m}}(p)=\varepsilon_{m} (2\pi)^{\frac{d}{2}}
\frac{\widehat{H_{m}}(p)\widehat{G_{m}}(p)}{|p|^{2}+|p|^{4}}, \quad 1\leq m\leq N.
$$
Then we obtain the expression for the norm as
\begin{equation}
\label{un}
\|u_{m}\|_{L^{2}({\mathbb R}^{d})}^{2}={(2\pi)}^{d}{\varepsilon}_{m}^{2}\int_{{\mathbb R}^{d}}
\frac{|\widehat{H_{m}}(p)|^{2}|\widehat{G_{m}}(p)|^{2}}
{[|p|^{2}+|p|^{4}]^{2}}dp.
\end{equation}
We use the analog of bound  (\ref{fub})  applied to functions $H_{m}$ and $G_{m}$ with
$R\in (0, +\infty)$ as
$$
{(2\pi)}^{d}{\varepsilon}_{m}^{2}\int_{{\mathbb R}^{d}}
\frac{|\widehat{H_{m}}(p)|^{2}|\widehat{G_{m}}(p)|^{2}}
{[|p|^{2}+|p|^{4}]^{2}}dp\leq
$$
$$
{(2\pi)}^{d}{\varepsilon}_{m}^{2}\Big[\int_{|p|\leq R}\frac{|\widehat{H_{m}}(p)|
^{2}|\widehat{G_{m}}(p)|^{2}}{|p|^{4}}dp+\int_{|p|>R}\frac{|\widehat{H_{m}}
(p)|^{2}|\widehat{G_{m}}(p)|^{2}}{|p|^{4}}dp \Big]\leq
$$
\begin{equation}
\label{ub1}
 {\varepsilon}_{m}^{2}\|H_{m}\|_{L^{1}({\mathbb R}^{d})}^{2}\Bigg\{
\frac{1}{{(2\pi)}^{d}}\|G_{m}\|_{L^{1}({\mathbb R}^{d})}^{2}|S^{d}|\frac{R^{d-4}}
{d-4}+
\frac{\|G_{m}\|_{L^{2}({\mathbb R}^{d})}^{2}}{R^{4}}\Bigg\}.
\end{equation}
Recall the norm definitions (\ref{u1N}) and (\ref{l2vn}). By means of the triangle
inequality along with the fact that $v\in B_{\rho}$, we easily derive
$$
\|u_{0}+v\|_{L^{2}({\mathbb R}^{d}, {\mathbb R}^{N})}\leq
\|u_{0}\|_{H^{4}({\mathbb R}^{d},{\mathbb R}^{N})}+1.
$$
Sobolev embedding (\ref{e}) yields
$$
|u_{0}+v|_{{\mathbb R}^{N}}\leq c_{e}(\|u_{0}\|_{H^{4}({\mathbb R}^{d}, {\mathbb R}^{N})}+1).
$$
Note that
$$
G_{m}(x)=\int_{0}^{1}\nabla g_{m}(t(u_{0}(x)+v(x))).(u_{0}(x)+v(x))dt, \quad
1\leq m\leq N.
$$
Here and further down the dot will stand for the scalar product of two vectors in ${\mathbb R}^{N}$.
Let us use the ball $I$ defined in (\ref{i}). Thus,
$$
|G_{m}(x)|\leq \hbox{sup}_{z\in I}|\nabla g_{m}(z)|_{{\mathbb R}^{N}}
|u_{0}(x)+v(x)|_{{\mathbb R}^{N}}\leq M|u_{0}(x)+v(x)|_{{\mathbb R}^{N}},
$$
such that
$$
\|G_{m}\|_{L^{2}({\mathbb R}^{d})}\leq M\|u_{0}+v\|_{L^{2}({\mathbb R}^{d},{\mathbb R}^{N})}
\leq M(\|u_{0}\|_{H^{4}({\mathbb R}^{d},{\mathbb R}^{N})}+1).
$$
Clearly, for $t\in [0,1]$ and $1\leq m,j\leq N$,
$$
\frac{\partial g_{m}}{\partial z_{j}}(t(u_{0}(x)+v(x)))=\int_{0}^{t}\nabla
\frac{\partial g_{m}}{\partial z_{j}}(\tau(u_{0}(x)+v(x))).(u_{0}(x)+v(x))d\tau.
$$
Hence,
$$
\Big|\frac{\partial g_{m}}{\partial z_{j}}(t(u_{0}(x)+v(x)))\Big|\leq
\hbox{sup}_{z\in I}\Big|\nabla \frac{\partial g_{m}}{\partial z_{j}}\Big|_
{{\mathbb R}^{N}}|u_{0}(x)+v(x)|_{{\mathbb R}^{N}}\leq
$$
$$
\sum_{n=1}^{N}\Big\|\frac{\partial^{2} g_{m}}{\partial z_{n}\partial z_{j}}
\Big\|_{C(I)}|u_{0}(x)+v(x)|_{{\mathbb R}^{N}}.
$$
This implies that
$$
|G_{m}(x)|\leq |u_{0}(x)+v(x)|_{{\mathbb R}^{N}}\sum_{j=1}^{N}
\sum_{n=1}^{N}\Big\|\frac{\partial^{2} g_{m}}
{\partial z_{n}\partial z_{j}}\Big\|_{C(I)}|u_{0,j}(x)+v_{j}(x)|\leq
$$
$$
M|u_{0}(x)+v(x)|_{{\mathbb R}^{N}}^{2}.
$$
Therefore,
\begin{equation}
\label{G1}
\|G_{m}\|_{L^{1}({\mathbb R}^{d})}\leq M\|u_{0}+v\|_{L^{2}({\mathbb R}^{d},
{\mathbb R}^{N})}^{2}\leq M(\|u_{0}\|_{H^{4}({\mathbb R}^{d}, {\mathbb R}^{N})}+1)^{2}.
\end{equation}
This allows us to derive the upper bound for the right side of
(\ref{ub1}) given by
$$
{\varepsilon}_{m}^{2}M^{2}\|H_{m}\|_{L^{1}({\mathbb R}^{d})}^{2}
(\|u_{0}\|_{H^{4}({\mathbb R}^{d}, {\mathbb R}^{N})}+1)^{2}
\Bigg\{\frac{(\|u_{0}\|_{H^{4}({\mathbb R}^{d}, {\mathbb R}^{N})}+
1)^{2}|S^{d}|R^{d-4}}{(2 \pi)^{d}(d-4)}
+\frac{1}{R^{4}}\Bigg\},
$$
where $R\in (0, +\infty)$. Let us recall Lemma 1.4 to minimize
the expression above, such that
$$
\|u_{m}\|_{L^{2}({\mathbb R}^{d})}^{2}\leq {\varepsilon_{m}}^{2}M^{2}\|H_{m}\|_
{L^{1}({\mathbb R}^{d})}^{2}\times
$$
$$
(\|u_{0}\|_{H^{4}({\mathbb R}^{d}, {\mathbb R}^{N})}+1)^
{2+\frac{8}{d}}
\Bigg(\frac{|S^{d}|}{4}\Bigg)^{\frac{4}{d}}\frac{d}
{(d-4)(2\pi)^{4}}, \quad 1\leq m\leq N.
$$
Then
\begin{equation}
\label{ul2ub}
 \|u\|_{L^{2}({\mathbb R}^{d},{\mathbb R}^{N})}^{2}\leq {\varepsilon}^{2}M^{2}H^{2}    
(\|u_{0}\|_{H^{4}({\mathbb R}^{d}, {\mathbb R}^{N})}+1)^
{2+\frac{8}{d}}
\Bigg(\frac{|S^{d}|}{4}\Bigg)^{\frac{4}{d}}
\frac{d}
{(d-4)(2{\pi})^{4}}.
\end{equation}
By means of (\ref{aux}),
$$
[-\Delta+{\Delta}^{2}]u_{m}(x)=\varepsilon_{m}
\int_{{\mathbb R}^{d}}H_{m}(x-y)G_{m}(y)dy, \quad 1\leq m\leq N, \quad 5\leq d\leq 7.
$$
We apply here the standard Fourier transform (\ref{f}),
the analog of bound (\ref{fub}) used for the function $G_{m}$
and estimate (\ref{G1}). This gives us
$$
\|{\Delta}^{2} u_{m}\|_{L^{2}({\mathbb R}^{d})}^{2}\leq \varepsilon_{m}^{2}
\|G_{m}\|_{L^{1}({\mathbb R}^{d})}^{2}
\|H_{m}\|_{L^{2}({\mathbb R}^{d})}^{2}\leq
$$
$$
\varepsilon^{2}M^{2}
(\|u_{0}\|_{H^{4}({\mathbb R}^{d}, {\mathbb R}^{N})}+1)^{4}
\|H_{m}\|_{L^{2}({\mathbb R}^{d})}^{2}.
$$
Thus,
\begin{equation}
\label{32}
\sum_{m=1}^{N}\|{\Delta}^{2} u_{m}\|_{L^{2}({\mathbb R}^{d})}^{2}\leq
\varepsilon^{2}M^{2}(\|u_{0}\|_{H^{4}({\mathbb R}^{d}, {\mathbb R}^{N})}+1)^{4}Q^{2}.
\end{equation}
Let us recall the definition of the norm (\ref{u1N}). By virtue of
(\ref{ul2ub}) and (\ref{32}), we have
$$
\|u\|_{H^{4}({\mathbb R}^{d}, {\mathbb R}^{N})}\leq
\varepsilon M(\|u_{0}\|_{H^{4}({\mathbb R}^{d}, {\mathbb R}^{N})}+1)^{2}\times
$$
\begin{equation}
\label{rh}
\Bigg[H^{2}(\|u_{0}\|_{H^{4}({\mathbb R}^{d}, {\mathbb R}^{N})}+1)^{\frac{8}{d}-2}
\Bigg(\frac{|S^{d}|}{4}\Bigg)^{\frac{4}{d}}
\frac{d}
{(d-4)(2 \pi)^{4}}
+Q^{2}\Bigg]
^{\frac{1}{2}}.
\end{equation}
Note that $\|u\|_{H^{4}({\mathbb R}^{d}, {\mathbb R}^{N})}\leq \rho$ for all the values of $\varepsilon$, which satisfy condition (\ref{eps}), 
so that $u\in B_{\rho}$ as well.

We suppose that for a certain $v\in B_{\rho}$ system (\ref{aux}) has
two solutions $u_{1,2}\in B_{\rho}$. Then their
difference $w(x):=u_{1}(x)-u_{2}(x) \in H^{4}({\mathbb R}^{d}, {\mathbb R}^{N})$.
Obviously, it satisfies the homogeneous system of equations
$$
[-\Delta+{\Delta}^{2}]w_{m}(x)=0, \quad 1\leq m\leq N.
$$
Evidently, the operator $l: H^{4}({\mathbb R}^{d})\to L^{2}({\mathbb R}^{d})$
introduced in (\ref{l}) does not possess any nontrivial zero modes. This means that
$w(x)\equiv 0$ in ${\mathbb R}^{d}$. Therefore, problem
(\ref{aux}) defines a map
$T_{g}: B_{\rho}\to B_{\rho}$ for all $\varepsilon$ satisfying assumption
(\ref{eps}).

Let us establish that this map is a strict contraction. For that purpose, we choose
arbitrarily $v_{1}, v_{2}\in B_{\rho}$. By virtue of the argument above,
$u_{1,2}:=T_{g}v_{1,2}\in B_{\rho}$ as well when $\varepsilon$ satisfies condition (\ref{eps}).
By means of (\ref{aux}), we have
\begin{equation}
\label{aux1}
[-\Delta+{\Delta}^{2}]u_{1,m}(x)=\varepsilon_{m} \int_{{\mathbb R}^{d}}
H_{m}(x-y)g_{m}(u_{0}(y)+v_{1}(y))dy,
\end{equation}
\begin{equation}
\label{aux2}
[-\Delta+{\Delta}^{2}]u_{2,m}(x)=\varepsilon_{m} \int_{{\mathbb R}^{d}}
H_{m}(x-y)g_{m}(u_{0}(y)+v_{2}(y))dy,
\end{equation}
where $1\leq m\leq N, \quad 5\leq d\leq 7$.
 Let us introduce
$$
G_{1, m}(x):=g_{m}(u_{0}(x)+v_{1}(x)), \quad G_{2, m}(x):=g_{m}(u_{0}(x)+v_{2}(x)),
\quad 1\leq m\leq N.
$$
We apply the standard Fourier transform (\ref{f}) to both sides of
systems  (\ref{aux1}) and (\ref{aux2}). This gives us
$$
\widehat{u_{1, m}}(p)=\varepsilon_{m} (2\pi)^{\frac{d}{2}}
\frac{\widehat{H_{m}}(p)\widehat{G_{1, m}}(p)}{|p|^{2}+|p|^{4}}, \quad
\widehat{u_{2, m}}(p)=\varepsilon_{m} (2\pi)^{\frac{d}{2}}
\frac{\widehat{H_{m}}(p)\widehat{G_{2, m}}(p)}{|p|^{2}+|p|^{4}}.
$$
Clearly,
\begin{equation}
\label{u12mn}
\|u_{1, m}-u_{2, m}\|_{L^{2}({\mathbb R}^{d})}^{2}=\varepsilon_{m}^{2}{(2\pi)}^{d}
\int_{{\mathbb R}^{d}} \frac{|\widehat{H_{m}}(p)|^{2}
|{\widehat{G_{1, m}}(p)}-{\widehat{G_{2, m}}}(p)|^{2}}
{[|p|^{2}+|p|^{4}]^{2}}dp.
\end{equation}
Let us use bound (\ref{fub}) to estimate the right side of (\ref{u12mn})  from above as
$$
\varepsilon_{m}^{2}{(2\pi)}^{d}\Bigg[\int_{|p|\leq R}\frac{|\widehat{H_{m}}(p)|^{2}
|{\widehat{G_{1, m}}(p)}-{\widehat{G_{2, m}}}(p)|^{2}}{|p|^{4}}dp+
$$
$$
\int_{|p|>R}\frac{|\widehat{H_{m}}(p)|^{2}
|{\widehat{G_{1, m}}(p)}-{\widehat{G_{2, m}}}(p)|^{2}}{|p|^{4}}dp \Bigg]\leq
$$
$$
\varepsilon^{2} \|H_{m}\|_{L^{1}({\mathbb R}^{d})}^{2} 
\Bigg\{\frac {\|G_{1, m}-G_{2, m}\|_{L^{1}({\mathbb R}^{d})}^{2}}{(2\pi)^{d}}|S^{d}|
\frac{R^{d-4}}{d-4}+
\frac{\|G_{1, m}-G_{2, m}\|_{L^{2}({\mathbb R}^{d})}^{2}}{R^{4}}\Bigg\},
$$
where $R\in (0,+\infty)$. Evidently, for
$1\leq m\leq N$, we can express
$$
G_{1, m}(x)-G_{2, m}(x)=\int_{0}^{1}\nabla g_{m}(u_{0}(x)+tv_{1}(x)+(1-t)v_{2}(x)).
(v_{1}(x)-v_{2}(x))dt.
$$
Note that for $t\in [0,1]$
$$
\|v_{2}+t(v_{1}-v_{2})\|_{H^{4}({\mathbb R}^{d}, {\mathbb R}^{N})}\leq
t\|v_{1}\|_{H^{4}({\mathbb R}^{d}, {\mathbb R}^{N})}+
$$
$$
(1-t)\|v_{2}\|_{H^{4}({\mathbb R}^{d}, {\mathbb R}^{N})}
\leq \rho.
$$
This means that $v_{2}+t(v_{1}-v_{2})\in B_{\rho}$. We easily obtain the inequality
$$
|G_{1, m}(x)-G_{2, m}(x)|\leq \hbox{sup}_{z\in I}|\nabla g_{m}(z)|_{{\mathbb R}^{N}}
|v_{1}(x)-v_{2}(x)|_{{\mathbb R}^{N}}
\leq M|v_{1}(x)-v_{2}(x)|_{{\mathbb R}^{N}}.
$$
Then
$$
\|G_{1, m}-G_{2, m}\|_{L^{2}({\mathbb R}^{d})}\leq M\|v_{1}-v_{2}\|_
{L^{2}({\mathbb R}^{d}, {\mathbb R}^{N})}\leq M\|v_{1}-v_{2}\|_
{H^{4}({\mathbb R}^{d}, {\mathbb R}^{N})}.
$$
Let us write
$\displaystyle{\frac{\partial g_{m}}{\partial z_{j}}(u_{0}(x)+tv_{1}(x)+(1-t)
v_{2}(x))}$  for $1\leq m,j\leq N$ as
$$
\int_{0}^{1}\nabla \frac{\partial g_{m}}{\partial z_{j}}
(\tau[u_{0}(x)+tv_{1}(x)+(1-t)v_{2}(x)]).[u_{0}(x)+tv_{1}(x)+(1-t)v_{2}(x)]d\tau.
$$
For $t\in [0,1]$, we derive
$$
\Big|\frac{\partial g_{m}}{\partial z_{j}}(u_{0}(x)+tv_{1}(x)+(1-t)v_{2}(x))\Big|\
\leq
$$
$$
\leq\sum_{n=1}^{N}\Bigg\|\frac{\partial^{2}g_{m}}{\partial z_{n}\partial z_{j}}
\Bigg\|_{C(I)}(|u_{0}(x)|_{{\mathbb R}^{N}}+t|v_{1}(x)|_{{\mathbb R}^{N}}+
(1-t)|v_{2}(x)|_{{\mathbb R}^{N}}).
$$
Therefore,
$$
|G_{1, m}(x)-G_{2, m}(x)|\leq
M|v_{1}(x)-v_{2}(x)|_{{\mathbb R}^{N}}
\Big(|u_{0}(x)|_{{\mathbb R}^{N}}+\frac{1}{2}|v_{1}(x)|_{{\mathbb R}^{N}}+
\frac{1}{2}|v_{2}(x)|_{{\mathbb R}^{N}}\Big).
$$
By virtue of the Schwarz inequality, we arrive at the upper bound on
the norm $\|G_{1, m}-G_{2, m}\|_{L^{1}({\mathbb R}^{d})}$, which is given by
$$
M\|v_{1}-v_{2}\|_{L^{2}({\mathbb R}^{d},{\mathbb R}^{N})}\Big(\|u_{0}\|_
{L^{2}({\mathbb R}^{d},{\mathbb R}^{N})}+\frac{1}{2}\|v_{1}\|_
{L^{2}({\mathbb R}^{d},{\mathbb R}^{N})}
+\frac{1}{2}\|v_{2}\|_{L^{2}({\mathbb R}^{d},{\mathbb R}^{N})}\Big)\leq
$$
\begin{equation}
\label{g12}
M\|v_{1}-v_{2}\|_{H^{4}({\mathbb R}^{d},{\mathbb R}^{N})}
(\|u_{0}\|_{H^{4}({\mathbb R}^{d}, {\mathbb R}^{N})}+1).
\end{equation}
Hence, the estimate from above on the norm
$\|u_{1, m}-u_{2, m}\|_{L^{2}({\mathbb R}^{d})}^{2}$ is equal to
$$
\varepsilon^{2}\|H_{m}\|_{L^{1}({\mathbb R}^{d})}^{2}
M^{2}\|v_{1}-v_{2}\|_{H^{4}({\mathbb R}^{d}, {\mathbb R}^{N})}^{2}
\Big\{\frac{(\|u_{0}\|_{H^{4}({\mathbb R}^{d}, {\mathbb R}^{N})}+1)^{2}
|S^{d}| R^{d-4}}{(2\pi)^{d}(d-4)}+\frac{1}{R^{4}}\Big\}.
$$
Let us minimize the expression above over $R\in (0,+\infty)$ using
Lemma 1.4. Thus, for $1\leq m\leq N$
$$
\|u_{1, m}-u_{2, m}\|_{L^{2}({\mathbb R}^{d})}^{2}\leq
\varepsilon^{2}\|H_{m}\|_{L^{1}({\mathbb R}^{d})}^{2}M^{2}
\|v_{1}-v_{2}\|_{H^{4}({\mathbb R}^{d}, {\mathbb R}^{N})}^{2}\times
$$
$$
(\|u_{0}\|_{H^{4}({\mathbb R}^{d}, {\mathbb R}^{N})}+1)^{\frac{8}{d}}
\Bigg(\frac{|S^{d}|}{4}\Bigg)^{\frac{4}{d}}
\frac{d}{(2 \pi)^{4}(d-4)}.
$$
Consequently,
$$
\|u_{1}-u_{2}\|_{L^{2}({\mathbb R}^{d}, {\mathbb R}^{N})}^{2}\leq \varepsilon^{2}H^{2}
 M^{2}\|v_{1}-v_{2}\|_{H^{4}({\mathbb R}^{d}, {\mathbb R}^{N})}^{2}\times
$$
\begin{equation}
\label{u12n}
(\|u_{0}\|_{H^{4}({\mathbb R}^{d}, {\mathbb R}^{N})}+1)^{\frac{8}{d}}
\Bigg(\frac{|S^{d}|}{4}\Bigg)^{\frac{4}{d}}
\frac{d}{(2{\pi})^{4}(d-4)}.
\end{equation}
Using (\ref{aux1}) and (\ref{aux2}) with $1\leq m\leq N$, we derive
$$
[-\Delta+{\Delta}^{2}](u_{1, m}(x)-u_{2, m}(x))
=\varepsilon_{m}
\int_{{\mathbb R}^{d}}H_{m}(x-y)[G_{1, m}(y)-G_{2, m}(y)]dy.
$$
Apply the standard Fourier transform (\ref{f}) along with bounds
(\ref{fub}) and (\ref{g12}). Hence,
$$
\|{\Delta}^{2} (u_{1, m}-u_{2, m})\|_
{L^{2}({\mathbb R}^{d})}^{2}\leq
\varepsilon^{2}\|G_{1, m}-G_{2, m}\|_{L^{1}({\mathbb R}^{d})}^{2}
\|H_{m}\|_{L^{2}({\mathbb R}^{d})}^{2}\leq
$$
$$
\varepsilon^{2}M^{2}\|v_{1}-v_{2}\|_{H^{4}({\mathbb R}^{d}, {\mathbb R}^{N})}
^{2}(\|u_{0}\|_{H^{4}({\mathbb R}^{d}, {\mathbb R}^{N})}+1)^{2}
\|H_{m}\|_{L^{2}({\mathbb R}^{d})}^{2}.
$$
This yields that
$$
\sum_{m=1}^{N}\|{\Delta}^{2}(u_{1, m}-u_{2, m})\|_
{L^{2}({\mathbb R}^{d})}^{2}\leq
$$
\begin{equation}
\label{d12}
\varepsilon^{2}M^{2}\|v_{1}-v_{2}\|_{H^{4}({\mathbb R}^{d}, {\mathbb R}^{N})}^{2}
(\|u_{0}\|_{H^{4}({\mathbb R}^{d}, {\mathbb R}^{N})}+1)^{2}Q^{2}.
\end{equation}
Let us use inequalities (\ref{u12n}) and (\ref{d12}) to estimate the norm
$\|u_{1}-u_{2}\|_{H^{4}({\mathbb R}^{d}, {\mathbb R}^{N})}$ from above by
$\varepsilon M(\|u_{0}\|_{H^{4}({\mathbb R}^{d}, {\mathbb R}^{N})}+1)\times$
\begin{equation}
\label{contr}
\Bigg\{\frac{H^{2}(\|u_{0}\|
_{H^{4}({\mathbb R}^{d}, {\mathbb R}^{N})}+1)^{\frac{8}{2}-2}d}
{(d-4){(2 \pi)}^{4}}
\Bigg(\frac{|S^{d}|}{4}\Bigg)^{\frac{4}{d}}
+Q^{2}\Bigg\}^
{\frac{1}{2}}\|v_{1}-v_{2}\|_
{H^{4}({\mathbb R}^{d}, {\mathbb R}^{N})}.
\end{equation}
It can be trivially checked that for all the values of $\varepsilon$ satisfying condition
(\ref{eps}), the constant in the right side of (\ref{contr}) is less than one.
This means that the map $T_{g}: B_{\rho}\to B_{\rho}$ defined by the system of equations
(\ref{aux}) is a strict contraction.
Its unique fixed point $u_{p}$ is the only solution of system 
(\ref{pert}) in the ball $B_{\rho}$. The resulting
$u\in H^{4}({\mathbb R}^{d}, {\mathbb R}^{N})$ given by (\ref{r}) solves
problem (\ref{p}). Recall formula (\ref{rh}). Clearly, $u_{p}$ converges
to zero in the $H^{4}({\mathbb R}^{d}, {\mathbb R}^{N})$ norm as
$\varepsilon\to 0$.      \hfill\lanbox

\bigskip

Let us turn our attention to the proof of the second main proposition of our
article.

\bigskip


\setcounter{section}{3}
\setcounter{equation}{0}

\centerline{\bf 3. The continuity of the resulting solution}

\bigskip

\noindent
{\it Proof of Theorem 1.5.} Evidently, for all the values of
$\varepsilon$, which satisfy inequality (\ref{eps}), we have
$$
u_{p,1}=T_{g_{1}}u_{p,1}, \quad u_{p,2}=T_{g_{2}}u_{p,2}.
$$
Hence,
$$
u_{p,1}-u_{p,2}=T_{g_{1}}u_{p,1}-T_{g_{1}}u_{p,2}+T_{g_{1}}u_{p,2}-
T_{g_{2}}u_{p,2}.
$$
Obviously,
$$
\|u_{p,1}-u_{p,2}\|_{H^{4}({\mathbb R}^{d}, {\mathbb R}^{N})}\leq\|T_{g_{1}}u_{p,1}-T_{g_{1}}
u_{p,2}\|_{H^{4}({\mathbb R}^{d},  {\mathbb R}^{N})}+\|T_{g_{1}}u_{p,2}-T_{g_{2}}u_{p,2}\|_
{H^{4}({\mathbb R}^{d}, {\mathbb R}^{N})}.
$$
According to bound (\ref{contr}), we arrive at
$$
\|T_{g_{1}}u_{p,1}-T_{g_{1}}u_{p,2}\|_{H^{4}({\mathbb R}^{d}, {\mathbb R}^{N})}\leq \varepsilon \kappa
\|u_{p,1}-u_{p,2}\|_{H^{4}({\mathbb R}^{d}, {\mathbb R}^{N})},
$$
where $\kappa$ is introduced in (\ref{kap}). Note that
$\varepsilon \kappa<1$ because the map
$T_{g_{1}}: B_{\rho}\to  B_{\rho}$ is a strict contraction under our
conditions. Therefore,
\begin{equation}
\label{sigma}
(1-\varepsilon \kappa)\|u_{p,1}-u_{p,2}\|_{H^{4}({\mathbb R}^{d}, {\mathbb R}^{N})}\leq
\|T_{g_{1}}u_{p,2}-T_{g_{2}}u_{p,2}\|_{H^{4}({\mathbb R}^{d}, {\mathbb R}^{N})}.
\end{equation}
We use that for our fixed point $T_{g_{2}}u_{p,2}=u_{p,2}$ and define
$\eta(x):=T_{g_{1}}u_{p,2}(x)$. For $1\leq m\leq N, \ 5\leq d\leq 7$, we have
\begin{equation}
\label{12}
[-\Delta+{\Delta}^{2}]\eta_{m}(x)=\varepsilon_{m}
\int_{{\mathbb R}^{d}}
H_{m}(x-y)g_{1, m}(u_{0}(y)+u_{p,2}(y))dy,
\end{equation}
\begin{equation}
\label{22}
 [-\Delta+{\Delta}^{2}]u_{p,2,m}(x)=\varepsilon_{m}
\int_{{\mathbb R}^{d}}H_{m}(x-y)g_{2, m}(u_{0}(y)+u_{p,2}(y))dy.
\end{equation}
Let us define
$$
G_{1,2,m}(x):=g_{1,m}(u_{0}(x)+u_{p,2}(x)), \quad G_{2,2,m}(x):=g_{2,m}
(u_{0}(x)+u_{p,2}(x)).
$$
We apply the standard Fourier transform (\ref{f}) to both sides of the
systems of equations (\ref{12}) and (\ref{22}). This gives us
$$
\widehat{\eta_{m}}(p)=\varepsilon_{m} (2 \pi)^{\frac{d}{2}}
\frac{\widehat{H_{m}}(p)
\widehat{G_{1,2,m}}(p)}{|p|^{2}+|p|^{4}}, \quad
\widehat{u_{p,2,m}}(p)=\varepsilon_{m} (2 \pi)^{\frac{d}{2}}\frac{\widehat{H_{m}}
(p)\widehat{G_{2,2,m}}(p)}{|p|^{2}+|p|^{4}}.
$$
Therefore,
\begin{equation}
\label{xiu}
\|\eta_{m}-u_{p,2,m}\|_{L^{2}({\mathbb R}^{d})}^{2}=
\varepsilon_{m}^{2}{(2\pi)}^{d}
\int_{{\mathbb R}^{d}}\frac{|\widehat{H_{m}}(p)|^{2}|
\widehat{G_{1,2,m}}(p)-\widehat{G_{2,2,m}}(p)|^{2}}{[|p|^{2}+|p|^{4}]^{2}}
dp.
\end{equation}
Let us use (\ref{fub}) to obtain the estimate from above on the right side of (\ref{xiu}), so that
$$
\varepsilon_{m}^{2}{(2\pi)}^{d}\Bigg[\int_{|p|\leq R}\frac{|{\widehat{H_{m}}}(p)|^{2}|\widehat{G_{1,2,m}}(p)-\widehat{G_{2,2,m}}(p)|^{2}}{|p|^{4}}dp+
$$
$$
\int_{|p|>R}\frac{|{\widehat{H_{m}}}(p)|^{2}|\widehat{G_{1,2,m}}(p)-
\widehat{G_{2,2,m}}(p)|^{2}}{|p|^{4}}dp \Bigg]\leq
\varepsilon^{2}\|H_{m}\|_{L^{1}({\mathbb R}^{d})}^{2}\times
$$
$$
\Bigg\{\frac{|S^{d}|}{(2{\pi})^{d}}
\frac{\|G_{1,2,m}-G_{2,2,m}\|_{L^{1}({\mathbb R}^{d})}^{2}R^{d-4}}{d-4}+
\frac{\|G_{1,2,m}-G_{2,2,m}\|_{L^{2}({\mathbb R}^{d})}^{2}}{R^{4}}\Bigg\}
$$
with $R\in (0, +\infty)$. Note that
$$
G_{1,2,m}(x)-G_{2,2,m}(x)=\int_{0}^{1}\nabla[g_{1,m}-g_{2,m}](t(u_{0}(x)+u_{p,2}(x))).
(u_{0}(x)+u_{p,2}(x))dt.
$$
Hence,
$$
|G_{1,2,m}(x)-G_{2,2,m}(x)|\leq \|g_{1,m}-g_{2,m}\|_{C^{2}(I)}
|u_{0}(x)+u_{p,2}(x)|_{{\mathbb R}^{N}}.
$$
This yields
$$
\|G_{1,2,m}-G_{2,2,m}\|_{L^{2}({\mathbb R}^{d})}\leq \|g_{1,m}-g_{2,m}\|_{C^{2}(I)}
\|u_{0}+u_{p,2}\|_{L^{2}({\mathbb R}^{d}, {\mathbb R}^{N})}\leq
$$
$$
\|g_{1,m}-g_{2,m}\|_{C^{2}(I)}(\|u_{0}\|_{H^{4}({\mathbb R}^{d}, {\mathbb R}^{N})}+1).
$$
It can be easily verified that for $1\leq m, j\leq N$ and $t\in [0,1]$, we have
$$
\frac{\partial}{\partial z_{j}}(g_{1,m}-g_{2,m})(t(u_{0}(x)+u_{p,2}(x)))=
$$
$$
\int_{0}^{t}\nabla \Big[\frac{\partial}{\partial z_{j}}(g_{1,m}-g_{2,m})\Big]
(\tau(u_{0}(x)+u_{p,2}(x))).(u_{0}(x)+u_{p,2}(x))d\tau.
$$
Thus,
$$
\Big|\frac{\partial}{\partial z_{j}}(g_{1,m}-g_{2,m})(t(u_{0}(x)+u_{p,2}(x)))\Big|
\leq
$$
$$
 \sum_{n=1}^{N}\Bigg\|\frac{\partial^{2}(g_{1,m}-g_{2,m})}{\partial z_{n}
\partial z_{j}}\Bigg\|_{C(I)}|u_{0}(x)+u_{p,2}(x)|_{{\mathbb R}^{N}}.
$$
Evidently,
$$
|G_{1,2,m}(x)-G_{2,2,m}(x)|\leq \|g_{1,m}-g_{2,m}\|_{C^{2}(I)}
|u_{0}(x)+u_{p,2}(x)|_{{\mathbb R}^{N}}^{2},
$$
such that
$$
\|G_{1,2,m}-G_{2,2,m}\|_{L^{1}({\mathbb R}^{d})}\leq
\|g_{1,m}-g_{2,m}\|_{C^{2}(I)}\|u_{0}+u_{p,2}\|_{L^{2}({\mathbb R}^{d}, {\mathbb R}^{N})}^{2}\leq
$$
\begin{equation}
\label{G1222}
\|g_{1,m}-g_{2,m}\|_{C^{2}(I)}(\|u_{0}\|_{H^{4}({\mathbb R}^{d}, {\mathbb R}^{N})}+1)^{2}.
\end{equation}
Let us use the inequalities above to derive the upper bound on the norm
$\|\eta_{m}-u_{p,2,m}\|_{L^{2}({\mathbb R}^{d})}^{2}$. It is equal to
$\varepsilon^{2}\|H_{m}\|_{L^{1}({\mathbb R}^{d})}^{2}(\|u_{0}\|_
{H^{4}({\mathbb R}^{d}, {\mathbb R}^{N})}+1)^{2} \times$
$$
\|g_{1, m}-g_{2, m}\|_{C^{2}(I)}^{2}
\Bigg[(\|u_{0}\|_{H^{4}({\mathbb R}^{d}, {\mathbb R}^{N})}+1)^{2}
\frac{|S^{d}|R^{d-4}}{(2 \pi)^{d} (d-4)}+\frac{1}{R^{4}}\Bigg].
$$
We minimize this quantity over $R\in (0, +\infty)$ by virtue of
Lemma 1.4 and arrive at
$$
\|\eta_{m}-u_{p,2,m}\|_{L^{2}({\mathbb R}^{d})}^{2}\leq
$$
$$
\varepsilon^{2}\|H_{m}\|_{L^{1}({\mathbb R}^{d})}^{2}
(\|u_{0}\|_{H^{4}({\mathbb R}^{d}, {\mathbb R}^{N})}+1)^{2+\frac{8}{d}}
\|g_{1, m}-g_{2, m}\|_{C^{2}(I)}^{2}
\Bigg(\frac{|S^{d}|}{4}\Bigg)^{\frac{4}{d}}
\frac{d}{(2\pi)^{4}(d-4)}.
$$
Recall definitions (\ref{h2}), (\ref{l2vn}), (\ref{gN}). We obtain that
$$
\|\eta-u_{p,2}\|_{L^{2}({\mathbb R}^{d}, {\mathbb R}^{N})}^{2}\leq
$$
\begin{equation}
\label{eup2l2n}
\varepsilon^{2}H^{2}(\|u_{0}\|_{H^{4}({\mathbb R}^{d}, {\mathbb R}^{N})}+1)^{2+\frac{8}{d}}
\|g_{1}-g_{2}\|_{C^{2}(I, {\mathbb R}^{N})}^{2}\Bigg(\frac{|S^{d}|}{4}\Bigg)^{\frac{4}{d}}
\frac{d}{(2\pi)^{4}(d-4)}.
\end{equation}
According to (\ref{12}) and (\ref{22}) with $1\leq m\leq N, \ 5\leq d\leq 7$, we have
$$
[-\Delta+{\Delta}^{2}]\eta_{m}(x)=
\varepsilon_{m} \int_{{\mathbb R}^{d}}H_{m}(x-y)G_{1,2,m}(y)dy,
$$
$$
[-\Delta+{\Delta}^{2}]u_{p,2,m}(x)=
\varepsilon_{m}\int_{{\mathbb R}^{d}}H_{m}(x-y)G_{2,2,m}(y)dy,
$$
such that
$$
[-\Delta+{\Delta}^{2}][\eta_{m}(x)-u_{p,2,m}(x)]=\varepsilon_{m}\int_{{\mathbb R}^{d}}H_{m}(x-y)[G_{1,2,m}(y)-G_{2,2,m}(y)]dy.
$$
Let us apply here the standard Fourier transform (\ref{f}) along with estimates (\ref{fub})
and (\ref{G1222}). This implies that
$$
\|{\Delta}^{2} [\eta_{m}(x)-u_{p,2,m}(x)]\|_{L^{2}({\mathbb R}^{d})}^{2}\leq
\varepsilon^{2}\|G_{1,2,m}-G_{2,2,m}\|_{L^{1}({\mathbb R}^{d})}^{2}
\|H_{m}\|_{L^{2}({\mathbb R}^{d})}^{2}\leq
$$
$$
\varepsilon^{2}\|g_{1, m}-g_{2, m}\|_{C^{2}(I)}^{2}
(\|u_{0}\|_{H^{4}({\mathbb R}^{d}, {\mathbb R}^{N})}+1)^{4}
\|H_{m}\|_{L^{2}({\mathbb R}^{d})}^{2}.
$$
We use definitions (\ref{q2}) and (\ref{gN}) to derive the upper bound
\begin{equation}
\label{sd2eul2}
\sum_{m=1}^{N}\|{\Delta}^{2} (\eta_{m}-u_{p,2,m})
\|_{L^{2}({\mathbb R}^{d})}^{2}\leq
\varepsilon^{2}
\|g_{1}-g_{2}\|_{C^{2}(I,{\mathbb R}^{N})}^{2}
(\|u_{0}\|_{H^{4}({\mathbb R}^{d}, {\mathbb R}^{N})}+1)^{4}Q^{2}.
\end{equation}
By virtue of (\ref{u1N})  and (\ref{l2vn}) along with (\ref{eup2l2n}) and (\ref{sd2eul2}), we arrive at
$$
\|\eta-u_{p,2}\|_
{H^{4}({\mathbb R}^{d}, {\mathbb R}^{N})}\leq \varepsilon
\|g_{1}-g_{2}\|_{C^{2}(I, {\mathbb R}^{N})}\times
$$
$$
(\|u_{0}\|_{H^{4}({\mathbb R}^{d}, {\mathbb R}^{N})}+1)^{2}\Bigg[
\frac{H^{2}(\|u_{0}\|_{H^{4}({\mathbb R}^{d},{\mathbb R}^{N})}+1)^{\frac{8}{d}-2}d}
{(d-4)(2\pi)^{4}}
\Bigg(\frac{|S^{d}|}{4}\Bigg)^{\frac{4}{d}}
+Q^{2}\Bigg]^
\frac{1}{2}.
$$
Recall inequality (\ref{sigma}). The norm
$\|u_{p,1}-u_{p,2}\|_{H^{4}({\mathbb R}^{d}, {\mathbb R}^{N})}$ can be estimated from above
by
$$
\frac{\varepsilon}{1-\varepsilon \kappa}
(\|u_{0}\|_{H^{4}({\mathbb R}^{d}, {\mathbb R}^{N})}+1)^{2}\times
$$
$$
\Bigg[\frac{{H}^{2}
(\|u_{0}\|_{H^{4}({\mathbb R}^{d}, {\mathbb R}^{N})}+1)^{\frac{8}{d}-2}d}{(d-4)
(2 \pi)^{4}}
\Bigg(\frac{|S^{d}|}{4}\Bigg)^{\frac{4}{d}}+
Q^{2}\Bigg]^{\frac{1}{2}}\|g_{1}-g_{2}\|_{C^{2}(I, {\mathbb R}^{N})}.
$$
By means of formulas (\ref{kap}) and (\ref{jres}) we complete the proof of the
theorem. \hfill\lanbox

\bigskip


\setcounter{section}{4}
\setcounter{equation}{0}

\centerline{\bf 4. Auxiliary results}

\bigskip
\bigskip

Let us establish the solvability for the linear Poisson type equation involving the sum of the negative Laplacian and the bi-Laplacian 
in the left side and a square integrable right side 
\begin{equation}
\label{lp}
[-\Delta+{\Delta}^{2}]\phi(x)=f(x), \quad x\in {\mathbb R}^{d}, \quad d\geq 5.
\end{equation}
The technical statement below can be trivially proved by virtue of the standard Fourier transform (\ref{f}) applied to both sides of
(\ref{lp}).

\bigskip

\noindent
{\bf Lemma 4.1.} {\it  Suppose that $f(x): {\mathbb R}^{d}\to {\mathbb R}, \ d\geq 5$ is nontrivial and 
$f(x)\in L^{1} ({\mathbb R}^{d})\cap L^{2} ({\mathbb R}^{d})$.
Then problem (\ref{lp}) has a unique solution
$\phi(x)\in H^{4}({\mathbb R}^{d})$.}

\bigskip

\noindent
{\it Proof.} Let us demonstrate that if
$\phi(x)\in L^{2} ({\mathbb R}^{d})$ is a solution of equation (\ref{lp}) with
a square integrable right side, it will belong to
$H^{4} ({\mathbb R}^{d})$ as well. For that purpose, we apply the
standard Fourier transform (\ref{f}) to both sides of (\ref{lp}) and arrive at
$$
(|p|^{2}+|p|^{4})\widehat{\phi}(p)=\widehat{f}(p)\in
L^{2} ({\mathbb R}^{d}).
$$
Thus,
$$
\int_{{\mathbb R}^{d}}[|p|^{2}+|p|^{4}]^{2}|\widehat{\phi}(p)|^{2}dp<\infty.
$$
Evidently, the equality
$$
\|{\Delta}^{2}\phi\|_{L^{2} ({\mathbb R}^{d})}^{2}=\int_{{\mathbb R}^{d}}|p|^{8}
|\widehat{\phi}(p)|^{2}dp<\infty
$$
is valid.
Therefore, ${\Delta}^{2}\phi\in L^{2} ({\mathbb R}^{d})$. Recall
the definition of the norm (\ref{n}). We derive that
$\phi(x)\in H^{4} ({\mathbb R}^{d})$.

Let us establish  the uniqueness of solutions for our problem. Suppose that equation (\ref{lp}) possesses two solutions
$\phi_{1, 2}(x)\in H^{4} ({\mathbb R}^{d})$. Then the difference function
$w(x):=\phi_{1}(x)-\phi_{2}(x)\in H^{4} ({\mathbb R}^{d})$. Clearly, it solves the homogeneous problem
$$
[-\Delta+{\Delta}^{2}]w=0.
$$
Note that the operator $l: H^{4} ({\mathbb R}^{d})\to L^{2} ({\mathbb R}^{d})$ defined in (\ref{l})  does not have
any nontrivial zero modes. Therefore, $w(x)$ vanishes in
${\mathbb R}^{d}$, which gives the uniqueness of solutions for equation (\ref{lp}).

Let us apply the standard Fourier transform (\ref{f}) to both sides of
problem (\ref{lp}). This yields
\begin{equation}
\label{fihpfp}
\widehat{\phi}(p)=\frac{\widehat{f}(p)}{|p|^{2}+|p|^{4}}\chi_{\{|p|\leq 1\}}+
\frac{\widehat{f}(p)}{|p|^{2}+|p|^{4}}\chi_{\{|p|>1\}}.
\end{equation}
In formula (\ref{fihpfp}) and below $\chi_{A}$ will stand for the
characteristic function of a set $A\subseteq {\mathbb R}^{d}$.

Obviously, the second term in the right side of (\ref{fihpfp}) can be bounded
from above in the absolute value by
$\displaystyle{\frac{|\widehat{f}(p)|}{2}\in L^{2}({\mathbb R}^{d})}$ via
the given condition.

Recall inequality (\ref{fub}). Hence, the first term in the right side of (\ref{fihpfp})
can be easily estimated from above in the absolute value by
\begin{equation}
\label{fxl1p2s2}
\frac{\|f(x)\|_{L^{1}({\mathbb R}^{d})}}{(2\pi)^{\frac{d}{2}}|p|^{2}}
\chi_{\{|p|\leq 1\}}.
\end{equation}
Clearly,  (\ref{fxl1p2s2}) is contained in $L^{2} ({\mathbb R}^{d})$ if $d\geq 5$. \hfill\lanbox.

\bigskip


\section*{Acknowledgements}

Vitali Vougalter is grateful to Israel Michael Sigal for the partial support
by the NSERC grant NA 7901. Vitaly Volpert has been supported by the RUDN
University Strategic Academic Leadership Program.

\end{document}